\documentclass[]{amsart}
\usepackage{amssymb}
\usepackage{amsmath}
\usepackage{amsthm}
\usepackage{amscd}
\usepackage{color}
\usepackage{epsfig}
\usepackage{graphicx,psfrag}
\usepackage[all]{xy}

\definecolor{darkgreen}{cmyk}{1,0,1,.2}

\newtheorem{thm}{Theorem}[section]
\newtheorem{thm*}{Theorem*}
\newtheorem{proposition}[thm]{Proposition}
\newtheorem{corollary}[thm]{Corollary}

\newtheorem{lemma}[thm]{Lemma}

\theoremstyle{definition}
\newtheorem{defn}[thm]{Definition}
\newtheorem{exa}[thm]{Example}

\theoremstyle{remark}

\def\square{\hfill${\vcenter{\vbox{\hrule height.4pt \hbox{\vrule
width.4pt
height7pt \kern7pt \vrule width.4pt} \hrule height.4pt}}}$}

\def\R{{\mathbb R}}

\def\Z{{\mathbb Z}}
\def\N{{\mathbb N}}

\def\P{{\mathcal P}}

\def\A{{\mathcal A}}

\def\HH{{\mathcal H}}

\def\diam{{\rm{diam}}}

\def\con{{\rm{Cone}}_{\omega}}

\def\ulim{\lim_\omega}

\def\<{\langle}
\def\>{\rangle}
\def\lf{\lfloor}
\def\rf{\rfloor}

\long\def\Restate#1#2#3#4{
\medskip\par\noindent
{\bf #1 \ref{#2} #3}
{\it #4}\par\medskip
}

\DeclareMathOperator{\CAT}{CAT}
\newcommand\catzero{$\CAT(0)$ }

\newcommand\infinity\infty
\newcommand\wt\widetilde
\newcommand\inject\hookrightarrow
\newcommand\union\cup

\newcommand{\co}{\colon\thinspace}
\newcommand\join\Lambda
\newcommand\cross\times

\newcommand\lub\vee

\newcommand\glb\wedge

\renewcommand\paragraph[1]{\bigskip\noindent\textbf{#1} }

\def\RAAG{right-angled Artin group }

\def\G{\Gamma}
\def\AG{A_\Gamma}
\def\AJ{A_J}
\def\XG{X_\Gamma}
\def\Gone{{\Gamma_1}}
\def\Gtwo{{\Gamma_2}}

\def\coneX{X^{\omega}}
\def\coneXG{\XG^{\omega}}
\def\coneAG{\AG^{\omega}}

\def\coned{d_{\omega}}

\newcommand\seq[1]{\mbox{\boldmath$#1$}}

\begin{document}
\title{Divergence and quasimorphisms of right-angled Artin groups}
\author{Jason Behrstock}
\author{Ruth Charney}
\date{January 2010}

\thanks {R. Charney was partially supported by NSF grant DMS 0705396.}
\thanks{J. Behrstock was partially supported by NSF grant DMS 0812513.}

\begin{abstract} We give a group theoretic characterization of
geodesics with superlinear divergence in the Cayley graph of a
right-angled Artin group $\AG$ with connected defining graph.  We use
this to determine when two points in an asymptotic cone of $\AG$ are
separated by a cut-point.  As an application, we show that if $\G$
does not decompose as the join of two subgraphs, then $\AG$ has an
infinite-dimensional space of non-trivial quasimorphisms.  By the work
of Burger and Monod, this leads to a superrigidity theorem for
homomorphisms from lattices into right-angled Artin groups.
\end{abstract}

\maketitle



\section{Introduction}

The \emph{divergence of a geodesic}, 
$\gamma\co[-\infty,\infty]\to X$, in a metric space, can roughly be 
thought of as the growth rate of a function from $\N$ to $\R$, 
which for each $N\in\N$ 
gives the length of the shortest path in $X\setminus 
B_{N}(\gamma(0))$ from $\gamma(-N)$ to $\gamma(N)$, where 
$B_{N}(\gamma(0))$ is the open ball of radius $N$ about $\gamma(0)$.
We refer to the \emph{divergence of a finitely generated group} to 
mean the largest divergence over all geodesics in a Cayley graph of 
$G$.

The divergence function has proven to be a useful tool in studying 
the large scale geometry of groups. Gersten classified geometric 
$3$--manifolds by their divergence \cite{Gersten:divergence3mflds} 
which allows one to distinguish quasi-isometry classes of 
$3$--manifold groups containing hyperbolic pieces from graph 
manifold groups \cite{KapovichLeeb:3manifolds}.  In addition, divergence
functions are closely related to cut-points in the asymptotic cones
of a group. 
Interest in the existence of such cut-points 
arose from Dru\c{t}u--Osin--Sapir's result that a group 
is relatively hyperbolic with respect to a collection of subgroups 
${\mathcal H}$ if and 
only if every asymptotic cone has a collection of cut-points with the
property that the maximal subsets of the asymptotic cone not separated
by any one of these cut-points arise from asymptotic cones of the
subgroups ${\mathcal H}$ \cite{DrutuSapir:TreeGraded}.

On the other hand, cut-points in asymptotic cones also arise in groups
which are not relatively hyperbolic.  To prove that any point in an
asymptotic cone of a mapping class group is a cut-point, the first
author showed that axes of pseudo-Anosov elements in a mapping class
group have super-linear divergence.  This also implies that these
directions are quasi-geodesically stable, or equivalently, Morse
geodesics.  \cite{Behrstock:asymptotic}.  (Alternate proofs have since
been obtained by \cite{DrutuMozesSapir:Divergence} and
\cite{DuchinRafi}).  More recently, Dru\c{t}u--Mozes--Sapir showed in
\cite{DrutuMozesSapir:Divergence} that a group has superlinear
divergence if and only if its asymptotic cones contain cut-points, and
that this occurs if and only if the group contains Morse geodesics.

In this paper we discuss divergence in right-angled Artin groups.
Given a finite, simplicial graph $\Gamma$, the \emph{right-angled Artin group} $\AG$ is the 
finitely presented group with generators corresponding to vertices of 
$\Gamma$ and relators of the form $x^{-1}y^{-1}xy$ whenever the 
vertices $x$ and $y$ of $\Gamma$ are connected by an edge. 
Right-angled Artin groups form 
a rich family of groups interpolating between $\Z^{n}$, the group 
corresponding to the complete graph on $n$ vertices, and the 
free group $F_n$, corresponding to the graph with $n$ 
vertices and no edges. 

If $\G_1$ and $\G_2$ are two graphs, their \emph{join} is the graph obtained by
connecting every vertex of $\G_1$ to every vertex of $\G_2$ by an edge.
Subgraphs of $\G$ that decompose as joins
are central to understanding divergence of geodesics. 
We define a notion of 
\emph{join length} of a geodesic, which measures the number of
cosets of join subgroups the geodesic passes through
(see Section~\ref{sec:joins} for the precise definition) and we prove,

\Restate{Theorem}{divergence and join length}{(Divergence and join length).}
{Let  $\Gamma$ be a connected graph and let $\alpha$ be a 
bi-infinite geodesic in $\AG$.  
Then $\alpha$ has linear divergence if and only if the 
join length of $\alpha$ is finite.}

The proof uses the action of $\AG$ on a \catzero cube complex,
$\XG$, the universal cover of the Salvetti complex of $\AG$.  We show that the join length
of a geodesic  $\alpha$ determines the behavior of the walls in $\XG$ crossed by $\alpha$. 

From the divergence theorem, we obtain the following complete characterization of 
when two points in an asymptotic cone of a right-angled Artin group
can be separated by a cut-point.

\Restate{Theorem}{classificationpieces}{(Classification of pieces).}
{Let $\Gamma$ be a connected graph. 
Fix a pair of points $\seq a, \seq b\in\coneAG$. 
The following are equivalent.
\begin{enumerate}
\item No point of $\coneAG$ separates $\seq a$ from $\seq b$.

\item There exist points $\seq a',\seq b'$ arbitrarily close
  to $\seq a,\seq b$, respectively for which the 
  join length between $\seq a',\seq b'$ 
  is finite.
 \end{enumerate}}

In the terminology of \cite{DrutuSapir:TreeGraded}, cut-points in an
asymptotic cone give rise to a \emph{tree-grading} whose \emph{pieces} 
are the
maximal subsets that cannot be separated by a point.  The above result
gives a complete description of the pieces in $\coneAG$.  Since
right-angled Artin groups are not relatively hyperbolic
\cite{BehrstockDrutuMosher:thick}, these pieces do not arise by taking
asymptotic cones of subgroups of $\AG$ \cite{DrutuSapir:TreeGraded}.
 
In \cite{BehrstockDrutuMosher:thick},  Behrstock--Dru\c{t}u--Mosher 
introduce a notion of \emph{algebraic thickness} of a group.
Theorem~\ref{classificationpieces} 
shows that for a connected graph $\Gamma$, $\AG$ is algebraically
thick of order zero if $\Gamma$ is a join, and otherwise it is 
algebraically thick of
order at least one with respect to the set of maximal join subgroups. 
It was established in \cite[Corollary~10.8]{BehrstockDrutuMosher:thick}, 
that, except for $\Z$, right-angled Artin groups with connected 
presentation graph are thick of order at most one. Together, 
these two results show that if $\AG$ is a join, then it is 
algebraically thick of order exactly zero, and otherwise it is 
algebraically thick of order exactly one.

Our main application of divergence is to  show that 
subgroups of right-angled Artin groups have many non-trivial
quasi-morphisms. 
A function $\phi\co G \to \R$ is a \emph{homogeneous quasimorphism} 
if $\phi(g^{n})=n\phi(g)$ for all $n>0$, and there exists a constant $D\geq 0$ such that 
$$|\phi(gh)-\phi(g)-\phi(h)|\leq D$$
for every $g,h\in G$. 
The vector space of homogeneous quasimorphisms, modulo the subspace of
true homomorphisms,
is denoted $\widetilde{QH}(G)$ and is related to the 
bounded cohomology of $G$. 
Bestvina and Fujiwara~\cite{BestvinaFujiwara:symspacebddcohom} have shown that for
group actions on a \catzero space, satisfying a weak discontinuity
property, the existence of
rank-one isometries (i.e., hyperbolic isometries with an axis 
 not bounding a half-plane) gives rise to non-trivial
quasimorphisms.  Using their results we prove,

\Restate{Theorem}{quasimorphisms} {(Rank-one geodesics and 
Quasimorphisms).}  
{     If $G\subseteq\AG$ is any non-cyclic, finitely generated subgroup which 
   is not contained in a conjugate of a join subgroup, then $G$ 
   contains an element which acts as a rank-one isometry of $\XG$. 
   In this case,  $\widetilde{QH}(G)$ is infinite dimensional.
   }

Right-angled Artin groups have been shown to have an extremely rich 
family of subgroups, cf.\    \cite{BestvinaBrady:MorseTheory}, \cite{HaglundWise},
\cite{CrispWiest}. In contrast, using Theorem~\ref{quasimorphisms} and 
Burger-Monod's result on nonexistence of quasimorphisms 
on higher rank lattices 
\cite{BurgerMonod:BddCohomLattices,BurgerMonod:BddCohomRigidity,
Monod:ContBddCohomBook}, we deduce:
   
   \Restate{Corollary}{superrigidity}{(Superrigidity with RAAG image).}
   {
   Let $\Lambda$ be an irreducible lattice in a connected 
   semisimple Lie group with finite center, no compact factors, and 
   rank at least 2. Then for any right-angled Artin group  
   $A_{\Gamma}$, every homomorphism 
   $\rho\co \Lambda \to A_{\Gamma}$ is trivial.
   }

   To the best of our 
   knowledge this is the most general statement of superrigidity for right-angled Artin groups,
 although many special cases follow from other  methods.   
   For example, for lattices satisfying Kazhdan's Property~(T), 
   Niblo--Reeves~\cite{NibloReeves:CAT0actions}
   showed that every action on a finite dimensional \catzero cube complex has a global fixed point.
   Since $\AG$ acts freely on the cube complex $\XG$, 
   any homomorphism of such a lattice into $\AG$ must be trivial.
   For non-uniform lattices, superrigidity follows from the Margulis 
   Normal Subgroup Theorem 
   \cite{Margulis:DiscreteSubgroups, Zimmer:Book}, since 
   every solvable subgroup of $\AG$ is virtually abelian.
   Other special cases follow from the work of Shalom  \cite{Shalom:Rigiditycommensurators}, 
   Monod \cite{Monod:superrigidsplittings}, and 
   Gelander--Karlsson--Margulis \cite{GelanderKarlssonMargulis}.

From Theorem~\ref{quasimorphisms} we also deduce the following, which P. Dani 
informed us she has independently established in joint work with 
A. Abrams, N. Brady, M. Duchin, A. Thomas and R. Young.

   \Restate{Corollary}{quaddivergence}{(Quadratic divergence).}
   {
    Let $\Gamma$ be a connected graph. $\AG$ has linear divergence if 
    and only if $\Gamma$ is a join; otherwise its divergence is 
    quadratic.
    }

    \bigskip
    
The authors would like to thank Koji Fujiwara  and Yehuda Shalom 
for helpful conversations. Also, Behrstock would like to thank 
Brandeis University and Charney would like to thank the Forschungsinstitut
f\"ur Mathematik at ETH, Zurich for their hospitality during the development 
of this paper.


\section{Walls}\label{sec:walls}

Let $\G$ be a finite, simplicial graph with vertex set $V$.  The \emph{\RAAG} associated to $\G$ is the group $\AG$ with presentation
\[
\AG=\langle V \mid vw=wv \hbox{ if  $v$ and $w$ are connected by an edge in $\G$}\rangle.
\]
Associated to any \RAAG $\AG$ is a CAT(0) cubical complex $\XG$ with a free action of $\AG$.
In this section we describe $\XG$ and investigate  the structure of walls in this complex.

Let $T$ be a torus of
dimension $|V|$ with edges labelled by the elements of $V$.  Let
$\overline X_\G$ denote the subcomplex of $T$ consisting of all faces
whose edge labels span a complete subgraph in $\G$ (or equivalently,
mutually commute in $\AG$).  $\overline X_\G$ is called the
\emph{Salvetti complex} for $\AG$.  It is easy to verify that the
Salvetti complex has fundamental group $\AG$ and that the link of the
unique vertex is a flag complex.  It follows that its universal cover,
$\XG$, is a CAT(0) cube complex with a free, cocompact action of
$\AG$.

If $\G'$ is a full subgraph of $\G$, then the inclusion $\G' \to \G$ induces an injective homomorphism
$\A_{\G'} \to \AG$ and an embedding $X_{\G'} \to \XG$.  This embedding
is locally geodesic, and hence (since $\XG$ CAT(0)) it is globally geodesic.  
We may thus view $X_{\G'}$ as a convex subspace of $\XG$.

We now recall some basic facts about walls in a CAT(0) cube
complex.  A \emph{wall} (or \emph{hyperplane}) in a CAT(0) cube
complex, $X$, is an equivalence
class of midplanes of cubes where the equivalence relation is
generated by the rule that two midplanes are related if they share a
face.  Each wall is a geodesic subspace of $X$ and separates $X$ into
two components.  Moreover, if a wall contains a (positive length)
segment of a geodesic $\gamma$, then it contains the entire geodesic
$\gamma$.

In the discussion that follows, we are interested in the relation between non-intersecting pairs of walls. The following terminology will be convenient. 

\begin{defn} Two walls $H_1, H_2$ in a CAT(0) cube complex are \emph{strongly separated}
if $H_1 \cap H_2 = \emptyset$ and no wall intersects both $H_1$ and $H_2$.
\end{defn}

Consider this definition in the context of a \RAAG $\AG$ and the
cube complex $\XG$.  For example, suppose
$\Gamma$ consists of two disjoint edges, so  $\AG$
is the free product $\Z^2 \ast \Z^2$.  In this case, the Salvetti
complex is the wedge of two tori, and its universal cover $\XG$
consists of flats which pairwise intersect in at most one vertex.
Since walls never contain vertices of $\XG$, they remain
entirely in one flat.  It follows that two walls are strongly
separated if and only if they lie in different flats.

At the other extreme, suppose $\G$ is a square, in which case 
$\AG= F_2 \times F_2$, the product of two free groups of rank 2, and
$\XG$ is the product of two trees $T_1 \times T_2$.  The walls
consist of trees of the form $T_1 \times \{y\}$ and $\{x\} \times T_2$
where $x$ and $y$ are midpoints of edges in $T_1$ and $T_2$
respectively.  It is now easy to see that no two walls are strongly
separated.

Now let $\AG$ be an arbitrary \RAAG and let $H_1$ and $H_2$ be two
walls in $\XG$.  Consider the set of all minimal length geodesics from
$H_1$ to $H_2$.  It follows from 
\cite[Proposition~II.2.2]{BridsonHaefliger} 
that the union of all such paths forms a
convex subspace of $\XG$ which we call the \emph{bridge} between
$H_1$ and $H_2$.

\begin{lemma}\label{bridge}  
If $H_1$ and $H_2$ are strongly separated, then the bridge $B$ between
them consists of a single geodesic from $H_1$ to $H_2$.
\end{lemma}

\begin{proof} It suffices to show that $B \cap H_1$ (and by symmetry $B \cap H_2$)  is a single point.
We first show that $B \cap H_1$ does not intersect any other wall $H$.
For suppose $x \in B \cap H_1 \cap H$.  Since $x
\in B$, it is the initial point of some minimal length geodesic
$\gamma$ from $H_1$ to $H_2$.  The initial segment of $\gamma$ lies in
some cube $\sigma$ of $\XG$ which contains midplanes in both $H_1$ and
$H$.  These midplanes span $\sigma$, hence the initial segment of
$\gamma$, which is orthogonal to $H_1$, must lie in $H$.  It follows
that all of $\gamma$ lies in $H$ and hence $H \cap H_2 \neq
\emptyset$.  This contradicts the assumption that $H_1$ and $H_2$ are
strongly separated.

Now every wall $H$ has an open neighborhood $N(H)$ isometric to $H \times (0,1)$,
namely the neighborhood consisting of the interiors of all cubes containing a midplane in $H$.
Then the same argument as above (using parallel copies of $H$ in $N(H)$)  
shows that $B \cap H_1 \cap N(H)$ must also be empty for all $H \neq H_1$.  
The only convex subsets of $H_1$ disjoint from every $N(H)$ are single vertices,
so it follows that $B \cap H_1$ is a single point.
 \end{proof}

 \begin{lemma}\label{distance}   There are universal constants $C, D >1$ (depending only on the dimension of $\XG$) such that  if $H_1$ and $H_2$ are strongly separated and $B$ is the bridge between them, then 
 \begin{enumerate}
 \item for any $x \in H_1$ and $y \in H_2$,
 $$ d(x,y) \geq  C^{-1}(d(x,B) + d(y,B)) -  d(H_1,H_2)-4$$
 \item  for any geodesic $\alpha$ in $\XG$, and any point $c$ on $\alpha$, if $H_1$ and $H_2$ intersect $\alpha$ inside the ball of radius $r$ about $c$, then the bridge $B$ is contained in the ball of radius $Dr$ about $c$.
 \end{enumerate}
 \end{lemma}
  
 \begin{proof}  
(1)  For any two points $x,y$ in $\XG$, define $d_H(x,y)$ to be the number of walls separating $x$ and $y$, or equivalently, the number of walls crossed by a geodesic from $x$ to $y$.  This distance function is quasi-isometric to the geodesic metric in $\XG$.  More precisely, 
 $$d(x,y)- C \leq d_H(x,y) \leq C d(x,y) + C $$
 where $C$ is the diameter of a maximal cube.

By Lemma~\ref{bridge}, $B$ consists of a single geodesic $\gamma$ from $H_1$ to $H_2$.  Let  $h_1, h_2$ be the endpoints of $\gamma$.  Let $\alpha$ be the geodesic  from $h_1$ to $x$,  and $\beta$ the geodesic from $y$ to $h_2$.  Note that $\alpha$ lies in $H_1$ and $\beta$ lies in $H_2$.  Since the strongly separated hypothesis guarantees that no wall crosses both $\alpha$ and $\beta$, the path $\alpha\gamma\beta$ crosses any given wall at most twice and $d_H(x,y)$ is the number of walls it crosses exactly once.  It follows that
 $$d_H(x,y)  \geq d_H(x,h_1) + d_H(y,h_2) - d_H(h_1,h_2). $$
Applying the inequalities above, we obtain
 \[
\begin{array}{rl}
d(x,y) &\geq C^{-1}d_H(x,y) -1 \\
&\geq  C^{-1}(d_H(x,h_1) + d_H(y,h_2) - d_H(h_1,h_2)) -1 \\
&\geq  C^{-1}(d(x,h_1) + d(y,h_2) - 2C) - d(h_1,h_2)-2 \\
& \geq  C^{-1}(d(x,B) + d(y,B)) -  d(H_1,H_2)-4.
\end{array}
\]

(2) Suppose $x=H_1 \cap \alpha$ and $y=H_2 \cap \alpha$ are in the ball of radius $r$ about $c$.  Then every point in $B$ is within 
$k=\frac{1}{2} (d(x,B) + d(y,B)) + d(H_1,H_2)$ of either $x$ or $y$ and hence within $k+r$ of $c$.
By part (1), $d(x,B)+d(y,B)$ is bounded by a linear function of $d(x,y)$, and by hypothesis, 
$d(H_1,H_2) \leq d(x,y) \leq 2r $ so $k$ is bounded by a linear function of $r$.  
\end{proof}

We now introduce the notion of \emph{divergence} for bi-infinite geodesics and discuss how the existence of strongly separated walls affects the divergence.  

\begin{defn}\label{defn divergence} Let $X$ be a geodesic metric
space.  Let $\alpha\co \R \to X$ be a bi-infinite geodesic in
$X$ and let $\rho$ be a linear function $\rho(r)= \delta r - \lambda$
with $0<\delta<1$ and $\lambda \geq 0$.  Define $div(\alpha,\rho)(r)$
= length of the shortest path from $\alpha(-r)$ to $\alpha(r)$ which
stays outside the ball of radius $\rho(r)$ about $\alpha(0)$ (or
$div(\alpha,\rho)(r)=\infty$ if no such path exists).  We say $\alpha$
has \emph{linear divergence} if for some choice of $\rho$,
$div(\alpha,\rho)(r)$ is bounded by a linear function of $r$, and that
$\alpha$ has \emph{super-linear divergence} otherwise.
\end{defn}

It is not difficult to verify that the definition of linear divergence is independent of
the choice of basepoint $\alpha(0)$.  We leave this as an exercise for the reader.

\begin{thm}\label{walls implies divergence}  Let $\alpha$ be a bi-infinite  geodesic in $\XG$ and let 
$\HH$ be the sequence of walls crossed by $\alpha$.  Suppose $\HH$ contains an infinite  subsequence 
$\{H_0, H_1, \dots \}$ such that for all $i$, $H_i$ is strongly separated from $H_{i+1}$.  Then $H_i$ is strongly separated from $H_j$ for all $i \neq j$ and $\alpha$ has superlinear divergence.  
\end{thm}

\begin{proof}  Let $\alpha^+$ and $\alpha^{-}$ denote the limit points of $\alpha$ in $\partial \XG$.  Since $H_i$ is disjoint from $H_{i+1}$, the half spaces $H_i^+$ containing $\alpha^+$ form a directed set $H_0^+ \subset H_1^+ \subset \dots $. Hence no two of these walls intersect and if some wall $H$ intersects both $H_i$ and $H_j$, $i<j$, then it must cross $H_{i+1}$, contradicting the strong separation of $H_i$ from $H_{i+1}$.  It follows that $H_i$ and $H_j$ are strongly separated for any $i<j$.  

Let $r'=\rho(r)$ and consider the ball $B_{r'}$ of radius $r'$ about $\alpha(0)$.  
Let $x_i= H_i \cap \alpha$.  By Lemma~\ref{distance}, for any $n$, we can choose 
$r$ large enough so that $B_{r'/2}$ contains $x_i$ for all $i \leq n$, as well as the  bridge between $H_{i-1}$ and $H_{i}$.  Let $\beta$ be any path from $\alpha(-r)$ to $\alpha(r)$ which 
stays outside $B_{r'}$.  Then $\beta$ must cross $H_0, H_1, \dots , H_{n}$ in a sequence of points $y_0, y_1, \dots y_{n}$.  Note that each $y_i$ is distance at least $r'/2$ from the bridges to the adjacent walls, hence by Lemma \ref{distance}, there is a universal constant $C$ such that 
 \[
\begin{array}{rl}
d(y_{i-1},y_i)  &\geq \frac{r'}{C}- (d(H_{i-1},H_i) -4 \\
&\geq \frac{r'}{C} - d(x_{i-1},x_i)-4
\end{array}
\]
  It follows that the length of $\beta$ satisfies
 \[
\begin{array}{rl}
|\beta| &\geq \sum d(y_{i-1},y_i)\\
& \geq \frac{nr'}{C}- 4n - d(x_0, x_n)\\
& \geq  \frac{nr'}{C}- 4n- r'
\end{array}
\]
Since $n \to \infty$ as $r \to \infty$, this proves the superlinear divergence of $\alpha$.
\end{proof}

The following example shows that the converse of the above theorem does not 
hold in complete generality.  However, when $\G$ is a connected graph, 
we will give a complete characterization of geodesics with superlinear divergence in
Theorem~\ref{divergence and join length} below.

\begin{exa} Suppose $\Gamma$ is disconnected, then $\AG$ splits as a free 
    product and $\XG$ splits as a wedge of spaces. Take any point 
    $p\in\XG$ whose removal disconnects $\XG$, and any pair of 
    geodesic rays $\gamma_{1}$ and  $\gamma_{2}$ emanating from $p$ 
    for which $\gamma_{1}\setminus\{p\}$ and 
    $\gamma_{2}\setminus\{p\}$ are in distinct components of 
    $\XG\setminus\{p\}$. Then the union of $\gamma_{1}$ and 
    $\gamma_{2}$ is a bi-infinite geodesic with super-linear 
    divergence (indeed infinite divergence, since $\gamma_{1}$ and 
    $\gamma_{2}$ can not be connected in the complement 
    of any ball around $p$). If we choose each of the $\gamma_{i}$ to 
    be contained in (the cube neighborhood of) a wall, then 
    $\alpha=\gamma_{1}\cup\gamma_{2}$ is 
    a geodesic with superlinear divergence such that no three walls
    crossed by $\alpha$ are pairwise strongly separated. 
\end{exa}


\section{Joins}\label{sec:joins}

In this section we give a group-theoretic interpretation of Theorem \ref{walls implies divergence}.
Choosing a vertex $x_0$ in $\XG$ as a basepoint, we may identify the
1--skeleton of $\XG$ with the Cayley graph of $\AG$ so that vertices
are labeled by elements of $\AG$ and edges by elements of the standard
generating set (namely the vertex set of $\Gamma$).  For a
generator $v$, let $e_v$ denote the edge from the basepoint $1$ to the vertex
$v$.  Any edge in $\XG$ determines a unique wall, namely the wall
containing the midpoint of that edge.  Denote by $H_v$ the wall
containing the midpoint of $e_v$.

For a cube in $\XG$, all of the parallel edges are labelled by the
same generator $v$.  It follows that all of the edges crossing a wall
$H$ have the same label $v$, and we call this a wall of \emph{type $v$}.
Since $\AG$ acts transitively on edges labeled $v$, a wall is of type
$v$ if and only if it is a translate of the standard wall $H_v$.

Let $lk(v)$ denote the subgraph of $\G$ spanned by the vertices adjacent to $v$ 
and $st(v)$ the subgraph spanned by $v$ and $lk(v)$.
The stabilizer of the wall $H_v$ is the group generated by
$lk(v)$, which we denote by $L_v$.  To see this, note that in any cube
containing the edge $e_v$, all other edges labeled $v$ are of the form
$ge_v$ for some $g \in L_v$.  An induction on the number of
cubes between $e_v$ and $e$ now shows that the same holds for any
edge $e$ which crosses $H_v$.

In what follows, for two subgroups $K$ and $L$ of $\AG$, we will use the notation $KL$ to mean the set of elements of $\AG$ which can be written as a product $kl$ for some $k \in K, l \in L$. In general, $KL$ is not a subgroup.

\begin{lemma} \label{ss_conditions}
 Let $H_1=g_1H_v$ and $H_2=g_2H_w$.  Then
\begin{enumerate}
\item $H_1$ intersects $H_2$ $\iff$  $v,w$ commute and $g_1^{-1}g_2 \in
L_v L_w$.
\item  $\exists ~H_3$ intersecting both $H_1$ and $H_2$ $\iff$   $\exists ~u \in st(v) \cap st(w)$ such that $g_1^{-1}g_2 \in L_v L_u  L_w$. 
\end{enumerate}
Here (2)  includes the case in which $H_3$ is equal to $H_1$ or $H_2$, hence
 $H_1$ and $H_2$ are strongly separated if and only if the conditions in (2) are not satisfied.
\end{lemma}

\begin{proof}  Without loss of generality, we may assume that $H_1=H_v$ and $H_2=gH_w$.

(1)  If $v,w$ commute, they span a cube in $\XG$, hence $H_v$ and $H_w$ intersect. 
Suppose $g=ab$, with $a \in L_v, b \in L_w$.  Then $H_v=a^{-1}H_v$ 
and $H_w=bH_w$, so translating by $a$, we see that $H_v$ intersects $gH_w$.  

Conversely, suppose $H_v$ intersects $gH_w$ in a cube $C$.  Then $C$
contains edges of type $v$ and of type $w$ hence $v$ and $w$ must
commute.  Moreover, $C$ is a translate $C=hC'$ of a cube $C'$ at the
basepoint containing the edges $e_v$ and $e_w$.  Since $e_v$ and
$he_v$ both intersect $H_v$, $h$ lies in $L_v$.  Since $ge_w$
and $he_w$ both intersect $gH_w$, $h^{-1}g$ lies in $L_w$.
Thus, $g \in L_v L_w$.

(2)  If $u \in st(v) \cap st(w)$ and $g=abc  \in L_v L_u  L_w$, then $H_v$ and $bH_w=bcH_w$ both intersect $H_u=bH_u$.  
Translating by $a$, we see that $H_v$ and $gH_w$ both intersect $aH_u$.

Conversely, suppose that $H_3=hH_u$ intersects both $H_1$ and $H_2$.
By part (1), $u$ must commute with both $v$ and $w$, so $u \in st(v)
\cap st(w)$.  Also by part (1), $h \in L_v L_u$ and
$h^{-1}g \in L_u L_w$, so $g \in L_v L_u L_w$.
\end{proof}

For two walls  $H_v$ and $gH_w$ to satisfy the conditions of (2),
both $w$ and the letters in $g$ must lie in a 2--neighborhood of $v$.
The converse is not true.  Consider for example the case of the Artin
group associated to a pentagon $\Gamma$ with vertices labeled (in
cyclic order) $a,b,c,d,e$.  Every vertex lies in a 2-neighborhood of
$a$, but we claim that $H_a$ and $daH_c$ are strongly separated.  This
follow from the lemma since $st(a) \cap st(c) = \{b\}$ and $da$ does not lie in 
$ L_aL_bL_c = \<e,b\> \< a,c\> \< b,d\>$.

To guarantee that no two walls in $\XG$ are strongly separated, we need a stronger condition.
Let $\Gone$ and $\Gtwo$ be (non-empty) graphs.  The \emph{join} of  $\Gone$ and $\Gtwo$ is the graph formed by joining every vertex of $\Gone$ to every vertex of $\Gtwo$ by an edge.  
The Artin group associated to such a graph splits as a direct product, $\AG = A_\Gone \times A_\Gtwo$
and $\XG$ splits as a metric product  $\XG= X_\Gone \times X_\Gtwo$.  The walls in $\XG$ are thus of the form $H_1 \times X_\Gtwo$ or $X_\Gone \times H_2$ for some wall $H_i$ in $X_{\Gamma_i}$.  
 Clearly, every wall of the first type intersects every wall of the second type, and it follows that no two walls are strongly separated.

Let $g \in \AG$ and let  $v_1v_2 \dots v_k$ be a minimal length word  representing $g$.  For $i < k$, set $g_i = v_1v_2 \dots v_i$.   Then the set of walls crossed by the edge path in $\XG$ from $x_0$ to $g x_0$ labelled $v_1v_2 \dots v_k$ is given by
$$ \HH=\{ H_{v_1}, g_1H_{v_2}, g_2H_{v_3} \dots g_{k-1}H_{v_k} \}.$$
A different choice of minimal word gives the same set of walls, changing only the order in which they are crossed.
Define the \emph{separation length} of $g$ to be
$$ \ell_S(g) = \max \{k \mid \textrm{$\HH$ contains $k$ walls which are pairwise strongly separated}\}.$$

 If $J$ is a complete subgraph of $\G$ which decomposes as a non-trivial join, then we call $A_J$ a \emph{join subgroup} of $\AG$. Define the \emph{join length} of $g$ to be
$$\ell_J(g) = \min \{ k \mid \textrm{$g=a_1 \dots a_k$ where $a_i$ lies in a join subgroup of $\AG$}\}. $$

If $\alpha$ is a (finite) geodesic in $\AG$, we can approximate $\alpha$ by an edge path which crosses the same set of walls as $\alpha$.  The word labeling this edge path determines an element $g_\alpha \in \AG$.  We define  $\ell_S(\alpha)= \ell_S(g_\alpha)$ and  $\ell_J(\alpha)= \ell_J(g_\alpha)$.
If $\alpha$ is a bi-infinite geodesic, and $\alpha_n$ denotes the restriction of 
$\alpha$ to the interval $[-n,n]$, we define the separation and join lengths of $\alpha$ to be
$$\ell_S(\alpha)= lim_{n \to \infty} \ell_S(\alpha_n) \quad\quad \ell_J(\alpha)= lim_{n \to \infty} \ell_J(\alpha_n).$$

\begin{lemma} \label{finite join length}
A bi-infinite geodesic $\alpha$ in $\AG$ has finite join length if and only if both the positive and negative rays of $\alpha$ eventually stay in a single join.  If every bi-infinite periodic geodesic in $\AG$ has finite join length, then $\G$ is itself a join. 
\end{lemma}

\begin{proof}   For any join $J$ in $\Gamma$,  $X_J$ is a convex subspace of $\XG$ so once $\alpha$ leaves $X_J$, it will never return, and similarly for translates of $X_J$.  If $\alpha$ has finite join length it  lies entirely in some finite set of these join subspaces and hence each ray must eventually remain in a single join.  The reverse implication is obvious.

For the second statement, suppose $\G$ is not a join.  Let $J$ be a
maximal join in $\G$ and let $v$ be a vertex not in $J$.  Let $g \in
A_J$ be the product of all the vertices in $J$ and consider the
bi-infinite geodesic $\alpha= \dots gvgvgvgv \dots$.  Note that no
vertex $w\in J \cup v$ commutes with both $J$ and $v$ since
otherwise, we would have $J\cup v$ contained in the join $st(w)$,
contradicting the maximality of $J$.
It follows that the tails of $\alpha$
must involve every vertex of $J \cup v$, hence by the first statement 
of the lemma, $\alpha$ has infinite join length.
 \end{proof}

In the proof of the previous lemma, we used the fact that for any vertex $v$ of $\G$,
$st(v)$ is always a join, namely it is the join of $\{v\}$ and $lk(v)$.  This fact plays a crucial role in the next lemma as well as in the proof of Theorem~\ref{divergence and join length} below.

\begin{lemma} \label{sep vs join length}
For any $g \in \AG$,
$$ \ell_S(g) \leq \ell_J(g) \leq 2\ell_S(g)+1. $$
Thus a bi-infinite geodesic has infinite join length if and only if it has infinite separation length.
\end{lemma}

\begin{proof} The first inequality follows from the observation above
that no two walls in a join are strongly separated.  For the second
inequality, fix a minimal word for $g$ and let $\HH$ be the sequence 
of walls crossed by the corresponding edge path as listed above.  
Set $H=H_{v_1}$ and let $H'=g_iH_{v_{i+1}}$ be the first
wall in the sequence strongly separated from $H$.  Then by
Lemma~\ref{ss_conditions}, $g_{i}$ lies in the product of three
link subgroups, $L_{v_1}L_{u_1} L_{v_{i+1}}$, for
some $u_1$, hence $g_{i+1}=g_iv_{i+1}$ lies in a product of
the three join groups generated by  $st(v_1)$, $st(u_1)$, and $st(v_{i+1}) $.  Now repeat
this argument starting with $H=g_iH_{v_{i+1}}$ and taking 
$H'=g_jH_{v_{j+1}}$ to be the next strongly separated wall (or
$H'=$ the last wall in $\HH$ if no more strongly separated walls
exist), to conclude that $g_{j+1}$ lies in the product of join subgroups
$$\< st(v_1)\> \<st(u_1) \> \< st(v_{i+1})\>\<st(u_2)\>\<st(v_{j+1})\>.$$
 Continuing this process, each new strongly separated wall adds two
 star subgroups.  Since we encounter at most $\ell_S(g)$ strongly
 separated walls, the inequality follows.
 \end{proof}


\section{The asymptotic cone}

The goal of this section is to understand the structure of the
asymptotic cones of $\AG$.  We begin by recalling some preliminaries
on asymptotic cones, tree graded spaces, and divergence; we refer the
reader to \cite{DriesWilkie},  \cite{DrutuSapir:TreeGraded}, and
\cite{DrutuMozesSapir:Divergence} for more details.

Let $(X,d)$ be a geodesic metric space.  Let $\omega$ be a non-principal
ultrafilter, $(o_n)$ a sequence of observation points in $X$, and
$(d_n)$ a sequence of scaling constants such that $\ulim d_n =
\infty$.  Then the asymptotic cone, $\con(X,(o_n),(d_n))$, is the
metric space consisting of equivalence classes of sequences $(a_n)$
satisfying $\ulim d(o_n,a_n)/d_n < \infty $, where two such sequences
$(a_n),(a'_n)$ represent the same point $\seq a$  
if and only if $\ulim d(a_n,a'_n)/d_n =0$, and the metric 
is given by $\coned (\seq a,\seq b)=\ulim d(a_n,b_n)/d_n$.

We will assume the observation points and scaling constants are fixed
and write $\coneX$ for $\con(X,(o_n),(d_n))$.  In general, 
$\coneX$ is a complete geodesic metric space. In the case where
$X$ has a cocompact group action, $\coneX$ is independent of choice of
observation points (but not, in general, of scaling constants) and is 
homogeneous.

A complete geodesic metric space is \emph{tree graded} if it contains 
a collection of closed subsets, $\P$, called \emph{pieces} such that 
the following three properties are satisfied: 
in each $P\in\P$, every pair of points is connected by a geodesic in 
$P$; any simple geodesic triangle is contained in some $P\in\P$; 
and each pair $P,P'\in\P$ is either disjoint or intersects in a single point. 
Dru\c{t}u--Osin--Sapir proved that a group is 
relatively hyperbolic if and only if all of its asymptotic cones are 
tree-graded with respect to pieces obtained by taking asymptotic cones 
of the peripheral
subgroups.  On the other hand, Behrstock--Dru\c{t}u--Mosher proved 
that right-angled Artin
groups are relatively hyperbolic if and only if their defining graph
is disconnected \cite{BehrstockDrutuMosher:thick}. 
In this section, we show that for connected defining
graphs, although the Artin group $\AG$ is not relatively hyperbolic,  
the asymptotic cones of $\AG$ still have a non-trivial
tree grading providing $\G$ is not a join. Moreover, although the 
pieces do not come from asymptotic cones of subgroups, 
they can be characterized group-theoretically (see
 Theorem~\ref{classificationpieces}). 

We begin by recalling the work of Drutu--Mozes--Sapir \cite{DrutuMozesSapir:Divergence} 
on divergence and cut-points. 

\begin{defn} Let $\rho(k)= \delta k - \lambda$ with $0<\delta<1$ and
$\lambda \geq 0$.  For points $a,b,c \in X$, set $k=d(c, \{a,b\}))$
and define $div(a,b,c;\rho)$ to be the length of the shortest path in
$X$ from $a$ to $b$ which lies outside the ball of radius $\rho(k)$
about $c$.  The \emph{divergence of $X$} with respect to $\rho$ is the
function
\[
Div(X,\rho)(r)=sup\{div(a,b,c;\rho) \mid a,b \in X, d(a,b) \leq r
\}.
\]
\end{defn}

For a biinfinite geodesic $\alpha$, the divergence function introduced in Section~\ref{sec:walls}
can be written as,
$$div(\alpha,\rho)(r)=div(\alpha(-r),\alpha(r),\alpha(0);\rho).$$
In particular, if $X$ has linear divergence, then every
bi-infinite geodesic in $X$ has linear divergence.

Drutu--Mozes--Sapir
establish the following correspondence between cut-points and divergence 
functions \cite[Lemma 3.14]{DrutuMozesSapir:Divergence}.

\begin{proposition}[\cite{DrutuMozesSapir:Divergence}]  \label{cut points}
 Let $\seq a=(a_n),  \seq b=(b_n)$, $\seq c=(c_n)$ be three points in $\coneX$,  
 and let $k = \coned (\seq c, \{\seq a, \seq b\})$.  Then
 $\seq c$ is a cut-point separating $\seq a$ from $\seq b$ if and only if for some $\rho$,
 \[
 \frac{lim_\omega div(a_n,b_n,c_n; \frac{\rho}{k})}{d_n }= \infty.
 \]
  \end{proposition}
  
  In particular, for a bi-infinite geodesic $\alpha$ in $X$, taking 
  $a_n=\alpha(-d_n), b_n=\alpha(d_n)$, and  $c_n=\alpha(0)$, 
 the proposition implies that  
 $\seq c$ is a cut-point separating $\seq a$ from $\seq b$ if and only if 
 $\alpha$ has super-linear divergence.

We say that $X$ is \emph{wide} if no asymptotic cone of $X$ has cut-points. 
In the case that $X$ is the Cayley graph of a group $G$, 
the proposition above leads to the following criterion for $G$ to be wide 
(see \cite[Proposition 1.1]{DrutuMozesSapir:Divergence}).  

\begin{proposition}[\cite{DrutuMozesSapir:Divergence}]  \label{wide}
A group $G$  is wide if and only if $Div(G,\rho)(r)$ is linear for $\rho(r)=\frac{1}{2} r -2.$  
\end{proposition}

In the case of $\AG$, the divergence of a bi-infinite geodesic is
determined by its join length.

\begin{thm}[Divergence and join length]\label{divergence and join length}
  Let  $\Gamma$ be a connected graph and let $\alpha$ be a bi-infinite geodesic in $\XG$.  Then $\alpha$ has linear divergence if and only if the join length of $\alpha$ is finite.  
  \end{thm}

\begin{proof}
If the join length of $\alpha$ is infinite, then by Lemma \ref{sep vs
join length}, so is the separation length.  By Theorem \ref{walls
implies divergence}, it follows that $\alpha$ has super-linear
divergence.

 Now suppose that the join length of $\alpha$ is finite.  We will show
 that $\alpha$ lies in a subspace of $\XG$ whose image in any
 asymptotic cone $\coneXG$ has no cut-points.  It then follows from
 the remarks following Proposition~\ref{cut points} that $\alpha$ has linear
 divergence.

By Lemma \ref{finite join length}, $\alpha$ lives entirely in a finite
union of join subspaces, that is, subspaces which are translates of
$X_J$ for some join $J$.  Since  $A_{J}$
decomposes as a direct product of infinite groups,  $X_J$ is wide.
Hence in any asymptotic cone
$\coneXG$, the cone on $gX_J$ gives rise to a subspace with no cut-points.  If $g'X_{J'}$ is another join subspace which shares a
geodesic line with $gX_J$, then the union of their asymptotic cones
contains a line in $\coneXG$ hence also has no cut-points.
 
Thus, it suffices to show that any two join subspaces $gX_J$
and $g'X_{J'}$ are connected by a sequence of join subspaces such that
consecutive subspaces share a line.  We will call this a
\emph{connecting sequence}.  By hypothesis, the graph $\G$ is
connected, so we can find a sequence of joins beginning at $J$ and
ending at $J'$ such that consecutive joins in the sequence share at least one
vertex in $\G$.  (For example, take a path from $J$ to $J'$ and take the sequence
of stars of the vertices along this path.) 
For $g=g'$,  it follows that there is a connecting
sequence from $gX_J$ to $gX_J'$.
  
For the general case, we may assume without loss of generality that
$g=1$.  Say $g'=a_1 \dots a_k$ where each $a_i$ lies in some join
$J_i$.  Then the observation above shows that there are connecting
sequences from $X_J$ to $X_{J_1}(=a_1X_{J_1}$), from $a_1X_{J_1}$ to
$a_1X_{J_2}(=a_1a_2X_{J_2}$), and so on to $g'X_{J_k}$, and finally,
from $g'X_{J_k}$, to $g'X_{J'}$.
\end{proof}

We now generalize the notion of join length to points in the
asymptotic cone $\coneAG$.  For two points $\seq a, \seq b\in\coneAG$,
we say that the pair $\seq a, \seq b$ (or the geodesic between
 them) has \emph{finite join length} if there exist representative sequences $(a_n),(b_n)$ for
which $\ulim \ell_J(a_n^{-1}b_n) < \infty$.

\begin{lemma} \label{cut points and join length}
Let $\seq a, \seq b$ be distinct points in $\coneAG$ and let $\seq\alpha$ be the geodesic between them.  If $\seq c$ is a point in the interior of $\seq\alpha$ which does not separate $\seq a$ from $\seq b$, then there exists a neighborhood of $\seq c$ in $\seq\alpha$ such that any two points in this neighborhood have finite join length. Moreover, the union of any two such neighborhoods of $\seq c$ also has this property. 
\end{lemma}

\begin{proof}   By Proposition~\ref{cut points}, the hypotheses imply that for some choice of linear function $\rho$, $div(a_n,b_n,c_n; \rho/k)$ is bounded by a linear function of $d_n$.  Let $D>0$ be as in Lemma \ref{distance} and let $k_n=d(a_n,b_n)$.  Since $\ulim \frac{k_n}{d_n} >0$,  we can choose  $\epsilon>0$ such that $D\epsilon < \frac{\rho(k_n)}{2d_n}$  for $\omega$-almost every $n$ .   Consider two sequences $(a'_n), (b'_n)$ lying within $\epsilon d_n$ of $(c_n)$.  For any two strongly separated walls between $a'_n$ and  $b'_n$, Lemma~\ref{distance} (2) implies that the bridge between them lies in the ball of radius $\frac{\rho(k_n)}{2}$ about $c_n$.  The number of such walls must be bounded (independent of $n$), since if not, then arguing as in Theorem~\ref{walls implies divergence},   we would deduce that $div(a_n,b_n,c_n; \rho/k)$ was super-linear.  It follows from Lemma~\ref{sep vs join length}  that the join length is also bounded.

To prove the last statement of the lemma, suppose $\seq c$ lies in 
two such intervals and let $\seq a', \seq b'$ be a point in each.  Then there exist two representative sequences $(c_n),(c'_n)$ for $\seq c$ such that the sequences $(a'^{-1}_nc_n)$ and $(c'^{-1}_nb'_n)$ have bounded join length.  If $(c^{-1}_nc'_n)$ has infinite join length, then $div(a'_n,b'_n,c_n;\rho/k')/d_n$ is unbounded for every $\rho$, so by Proposition~\ref{cut points}, $\seq c$ is a cut-point. This contradicts the hypotheses of the lemma, so we conclude that $(c^{-1}_nc'_n)$, and hence $(a'^{-1}_nb_n)$, has bounded join length.
\end{proof}

\begin{thm}[Classification of pieces]\label{classificationpieces}
   Let $\Gamma$ be a connected graph. 
   Fix a pair of points $\seq a, \seq b\in\coneAG$. 
   The following are equivalent.
   \begin{enumerate}
   \item No point of $\coneAG$ separates $\seq a$ from $\seq b$.

   \item There exist points $\seq a',\seq b'$ arbitrarily close
     to $\seq a,\seq b$, respectively for which the 
     join length between $\seq a',\seq b'$ 
     is finite.
    \end{enumerate}
\end{thm}

\begin{proof}
 Suppose (1) holds.  Let $\seq\alpha$ be the geodesic from $\seq a$ to $\seq b$.  Then by Lemma~\ref{cut points and join length}, every point on $\seq\alpha$ is contained in an open interval in which any two points have finite join length. Moreover, if two such intervals intersect, then their union also has this property. 
 It now follows easily that the maximal open interval of $\seq\alpha$ such that any two points have finite join length is the entire interior of $\seq\alpha$.

    Now suppose (2) holds.
    By hypothesis, for every $\epsilon>0$ there exist points 
    $\seq a'$ and $\seq b'$ with 
    representatives 
    $(a'_n)$ and $(b'_n)$ whose distances in 
    $\coneAG$ are less than $\epsilon$ from $\seq a$ and $\seq b$, 
    respectively, and the join length between 
    $\seq a'$ and $\seq b'$ is finite,  that is, there exists a constant $M$ 
    such that $\omega$-almost every
    $a'^{-1}_n b'_n$ is a product of at most
    $M$ subwords, each contained in a join subgroup.
    Hence the corresponding geodesic is contained in a finite
    sequence of join cosets.
    As in the proof of Theorem~\ref{divergence and join length} there 
    is a connecting sequence, that is,  
    a finite sequence of additional join cosets which we may add between any two of these
   to get an ordered sequence of join cosets where 
    each intersects the next in an infinite diameter set.  Starting 
    with $M$ join cosets, the construction in Theorem~\ref{divergence and join length} yields 
    a connecting sequence whose length is at most 
    $M\cdot\diam(\Gamma)$, where $\diam(\Gamma)$ is the diameter of 
    $\Gamma$.  Denote this sequence  by $\mathcal S_n$. 
    
   We will say that a coset of $A_J$ has join type $J$.  Since there
   are only finitely many joins in $\Gamma$, 
   the sequence of join types in $\mathcal S_n$ is the same 
   for $\omega$--almost every $n$.  Any two cosets of the same join type are 
   isometric, so for each $i$, the subspace of points $(x_n) \in \coneAG$ such that $x_n$ lies in the
   $i^{th}$ term of $S_n$ is isometric to $A^{\omega}_{J_i}$ for some join $J_i$.
   In particular, this subspace has no cut-points.
   Moreover, the intersection of any two consecutive subspaces in this sequence
   has infinite diameter.  It follows that their union,
   which contains $\seq a'$ and $\seq b'$, has no cut-points.
      
    Hence, $\seq a$ and $\seq b$ can be
    approximated arbitrarily closely by points 
    $\seq a'$ and $\seq b'$ which cannot be separated by a point.  
    Since not being separated by a point is a closed condition, 
    this completes the proof that (2) implies (1).
\end{proof}

As an immediate corollary we obtain,

\begin{corollary}  \label{wide equal join}
$\AG$ is wide if and only if $\G$ decomposes as a non-trivial join.
\end{corollary}

The following example shows that one cannot replace the second condition in the theorem by the simpler statement that the geodesic from $\seq a$ to $\seq b$ has finite join length.

\begin{exa} Suppose $x,y$ are two vertices in $\G$ that are not
contained in a join.  For simplicity, take the scaling constants for
$\coneAG$ to be $d_n=n$.  Let $a_n=1$ for all $n$ and let
$$b_n=x^{\lf \frac{n}{2}\rf}y^{\lf \frac{n}{4}\rf}x^{\lf
\frac{n}{8}\rf} y^{\lf \frac{n}{16}\rf }\dots $$
 Then the join length of from $\seq a$ to $\seq b$ is infinite.
 However, if we truncate each $b_n$ after $k$ terms, setting
$$b_n^{(k)}=x^{\lf \frac{n}{2}\rf}y^{\lf \frac{n}{4}\rf} \dots z^{\lf
\frac{n}{2^k}\rf}$$
where $z=x,y$ depending on whether $k$ is odd or even, we obtain a
point $\seq b^{(k)}$ in the asymptotic cone whose distance from
$\seq b$ is $\frac{1}{2^k}$ while the join length from $\seq a$ to
$\seq b^{(k)}$ is $k$.  It follows from the theorem above that $\seq
a$ and $\seq b$ lie in the same piece of the asymptotic cone, despite
the fact that the join length between them is infinite.
\end{exa}


\section{Divergence, quasimorphisms, and superrigidity}

Recall that a hyperbolic isometry of a proper \catzero space is called a 
\emph{rank-one isometry} if some axis of that isometry does not bound 
a half-plane. A \emph{quasimorphism} on a group $G$ is a function 
$\phi\co G \to \R$ for which there exists a constant $D(\phi)\geq 0$ 
such that  $$|\phi(gh)-\phi(g)-\phi(h)|\leq D(\phi)$$
for every $g,h\in G$. A quasimorphism is \emph{homogeneous} if for 
each $g\in G$ and each $n\in\N$, we have $\phi(g^{n})=n\phi(g)$. 
The set of homogeneous quasimorphisms on a given group $G$ form a 
vector space. The quotient of this vector space by homomorphisms from 
$G$ to $\R$ is denoted $\widetilde{QH}(G)$ and is isomorphic to the 
kernel of the map from the second bounded cohomology of $G$ (with $\R$ 
coefficients) to the second cohomology of $G$. (For details, see 
\cite{Gromov:Volume} and \cite{Calegari:scl}.)

Burger and Monod proved that there are no non-trivial 
homogeneous quasimorphisms on any irreducible lattice in a connected 
semisimple Lie group of rank at least 2 
with finite center and no compact factors 
\cite{BurgerMonod:BddCohomLattices,BurgerMonod:BddCohomRigidity, 
Monod:ContBddCohomBook}. On the 
other hand, several interesting 
families of groups, including non-elementary hyperbolic groups 
\cite{EpsteinFujiwara:BddcohomHypGps} and 
mapping class groups \cite{BestvinaFujiwara:boundedcohom}, 
have been shown to have infinite dimensional $\widetilde{QH}$. We now 
establish where right-angled Artin groups lay in this framework.

To prove the theorem we will need some basic facts about normal forms in right-angled Artin groups.
We refer the reader to \cite{Laurence} for details.  Let $V$ be the generating set for $\AG$ and let
$g$ be an element $\AG$.   A \emph{reduced word} for $g$ is a minimal length word in the free group $F(V)$ representing $g$.
Given an arbitrary word representing $g$, one can obtain a reduced word by a process of
``shuffling" (i.e. interchanging commuting elements) and canceling inverse pairs. Any two reduced
words for $g$ differ only by shuffling.  In particular, the \emph{support} of $g$, that is the set $Supp(g) \subseteq V$ of letters appearing in a reduced word for $g$, is independent of choice of reduced word.

An element of $\AG$ is called \emph{cyclically reduced} if it is of minimal length in its conjugacy class. 
For any $g \in \AG$, there exists a unique cyclically reduced element conjugate to $g$, which we denote by $\bar g$.  Given a reduced word $w$ representing $g$, we can find a reduced word $\bar w$ for $\bar g$ by shuffling $w$ to get a maximal length word $u$ such that $w=u\bar wu^{-1}$.  In particular, $g=\bar g$ if and only if every shuffle of $w$ is cyclically reduced as a word in the free group $F(V)$.

\begin{lemma} \label{centralizer}
 Let $g=\bar g$ be a cyclically reduced element of $\AG$.  Then the following are equivalent.
\begin{enumerate}
\item $g$ is contained in a join subgroup.
\item The centralizer of $g$ is non-cyclic.
\item The centralizer of $g$ is contained in a join subgroup.
\end{enumerate}
\end{lemma}

\begin{proof}
(1) implies (2) since the centralizer of any element $(g_1,g_2)$ of a
direct product $G_1 \times G_2$ is the product of the centralizers
$C_{G_1}(g_1) \times C_{G_2}(g_2)$.  (3) implies (1) is obvious, so it
remains only to prove that (2) implies (3).
 
For any subset $S \subset V$, let $lk(S)$ denote the (possibly empty)
set of vertices at distance 1 from every vertex of $S$.  It follows
from Servatius' Centralizer Theorem~\cite{Servatius} (see also Thm 1.2
of \cite{Laurence}) that the centralizer of a cyclically reduced
element $g$ lies in the subgroup generated by $Supp( g) \cup
lk(Supp(g))$ and that the centralizer is cyclic unless either
$lk(Supp(g))$ is nonempty or $Supp(g)$ decomposes as a join.  In
either case, $Supp( g) \cup lk(Supp(g))$ spans a join in $\G$.
\end{proof}

\begin{thm}[Rank-one geodesics and Quasimorphisms]  
    \label{quasimorphisms} 
   If $G\subseteq\AG$ is any non-cyclic, finitely generated subgroup which 
   is not contained in a conjugate of a join subgroup, then $G$ 
   contains an element which acts as a rank-one isometry of $\XG$. 
   In this case,  $\widetilde{QH}(G)$ is infinite dimensional.
\end{thm}

\begin{proof}   Let $g$ be an element of $\AG$ and let $\bar g$ denote the cyclic reduction of $g$. Then $\bar g^k$ is geodesic for all $k$.  If $\bar g$ does not lie in a join subgroup, then the axis for $\bar g$ has infinite join length, hence by Theorem~\ref{divergence and join length}, it has super-linear divergence.  It follows that the axis for $\bar g$ cannot bound a half-flat and the same holds for the axis of $g$ since it is a translate of the axis for $\bar g$.   Thus to prove the first statement of the theorem, it suffices to show that $G$ contains an element $g$ whose cyclic reduction $\bar g$ does not lie in a join subgroup.  

Choose an element $c \in G$ such that $Supp(\bar c)$ is 
maximal.  That is, if $g \in G$ 
has $Supp(\bar c)\subseteq Supp(\bar g)$, then $Supp(\bar c) = Supp(\bar g)$.   Conjugating $G$ if necessary,  we may assume without loss of generality that $c=\bar c$.  
If $c$ is not contained in a join,  we are done.  
So suppose $c$, and hence by Lemma~\ref{centralizer} the centralizer of $c$, is contained in 
a join subgroup $A_J$.  

By hypothesis, $G$ does not lie in a join, so there is some element $h \in G$ whose support is not contained in $J$. Consider an element of the form
 $x= c^khc^k \in G$.  We claim that for sufficiently large $k$, 
 $Supp(\bar x) \supsetneq Supp(c)$ contradicting the maximality assumption on $Supp(c)$.  To see this, note that since $c$ is cyclically reduced, cancellations can only occur between generators in $c$ and generators in $h$.  It follows that repeatedly multiplying $h$ on the left or right by $c$, can  result in at most finitely many cancellations and all canceled letters must lie in $Supp(c) \cap Supp(h)$.
Thus for $k$ sufficiently large, a reduced word for $x$ is of the form $c^i u c^j$ 
for some $i,j >0$ and some reduced word $u \in G$, and the support of $x$ satisfies
$Supp(x)=Supp(c^i u c^j) = Supp(c) \cup Supp(h)$.  

Now consider the cyclic reduction $\bar x$ of $x$.  
Write $u=a\bar u a^{-1}$, so $x=c^i a \bar u a^{-1} c^j$.  Since $c$ is cyclically reduced, the only
way $x$ can fail to be cyclically reduced is if some initial subword $a'$ of $a$ commutes with $c$. 
But in this case, we may conjugate $G$ by $a'^{-1}$ and repeat the argument replacing $h$ by 
$h'=a'^{-1}ha'$  to conclude that $x'=c^i\bar u c^j$ is cyclically reduced.  
Since $a'$ lies in the centralizer of $c$, it  
lies in $\AJ$,  whereas $h \notin \AJ$.  It follows that $h' \notin \AJ$, hence 
$Supp(\bar x')=Supp(x')  = Supp(c) \cup Supp(h') \supsetneq Supp(c)$, as claimed. 

For the second statement of the theorem, note that the 
action of any non-cyclic subgroup of $\AG$ on $\XG$ is
weakly properly discontinuous in the sense of
Bestvina--Fujiwara \cite{BestvinaFujiwara:boundedcohom}
 since it is a properly discontinuous action of a
non-virtually cyclic group on a \catzero space.  Since the action
of $G$ on $\XG$ contains a rank-one isometry, the Main Theorem
of \cite{BestvinaFujiwara:symspacebddcohom} implies that
$\widetilde{QH}(G)$ is infinite dimensional.
\end{proof}

\begin{corollary}[Superrigidity with RAAG image]\label{superrigidity}
    Let $\Lambda$ be an irreducible lattice in a connected 
    semisimple Lie group with finite center, no compact factors, and 
    rank at least 2. Then for any right-angled Artin group  
    $A_{\Gamma}$, every homomorphism 
    $\rho\co \Lambda \to A_{\Gamma}$ is trivial. 
\end{corollary}

\begin{proof} 
      We proceed by induction on the number $n$ of
    vertices in $\Gamma$.  If $n=1$, then $\AG$ is infinite cyclic and the 
    Margulis Normal Subgroup Theorem 
    \cite{Margulis:DiscreteSubgroups, Zimmer:Book}  
    implies that the image of $\rho$ is trivial.
    
    Suppose $n \geq 2$.   If the image
    $\rho(\Lambda)$ is not contained in the conjugate of a join subgroup, 
    then Theorem~\ref{quasimorphisms} yields a
    non-trivial quasimorphism on $\rho(\Lambda)$, and hence the 
    composition of this with $\rho$ gives a 
    nontrivial quasimorphism on $\Lambda$.  This contradicts
    Burger--Monod's result that such lattices do not admit any 
    non-trivial quasimorphisms 
    \cite{BurgerMonod:BddCohomLattices,BurgerMonod:BddCohomRigidity, 
    Monod:ContBddCohomBook}.
    
    Thus, up to conjugacy,  we may assume that $\rho(\Lambda)$ lies in a join 
    subgroup  $A_J=A_{\Gamma_{1}}\times A_{\Gamma_{2}}$. 
    Composing $\rho$ with the projections on each factor, gives two homomorphisms
    from $\Lambda$  to right-angled Artin groups with less than $n$ generators.  
    By induction, both of these homomorphisms are trivial, hence $\rho$ is also
    trivial.
\end{proof}

\begin{corollary}[Quadratic divergence]\label{quaddivergence}
    Let $\Gamma$ be a connected graph. $\AG$ has linear divergence if 
    and only if $\Gamma$ is a join; otherwise its divergence is 
    quadratic.
\end{corollary}
\begin{proof}
The statement that $\AG$ has linear divergence if and only if $\G$ is a join follows immediately
   from Proposition~\ref{wide} and Corollary~\ref{wide equal join}.
         
    For the second statement, consider the case that $\G$ is not a join. 
    Since $\XG$ is a locally-compact \catzero space, we can apply
    \cite[Proposition 3.3]{KapovichLeeb:3manifolds} which implies that
    any complete periodic geodesic in $\XG$ has divergence which is
    either linear or at least quadratic.  Theorem~\ref{quasimorphisms}
    states that $\AG$ contains a rank-one geodesic, hence the
    divergence of $\AG$ is at least quadratic.  
    
    On the other hand, in \cite{BehrstockDrutuMosher:thick} it is
    shown that for any connected graph $\G$, $\AG$ is algebraically
    thick of order at most 1.  It is not hard to show that any group
    which is algebraically thick of order 1 has at most quadratic
    divergence (a generalization of this fact for metric spaces which
    are thick of arbitrary degree will appear in
    \cite{BehrstockDrutu:thick2}).  We conclude that the divergence of
    $\AG$ is exactly quadratic.
\end{proof}


\end{document}